\documentclass[12pt]{amsart}
\usepackage{amssymb}
\usepackage{amsthm}
\usepackage{amsfonts}
\newtheorem{lemma}{\bf Lemma}[section]
\newtheorem{observation}[lemma]{\bf Observation}
\newcommand{\Pow}{{\mathcal P}}
\newcommand{\proofend}{{\hfill $\Box$}}
\newcommand{\card}{{\mathrm{card}}}

\begin{document}

\title{Counterexample to a conjecture of Aharoni and Korman}

\author{Dominic van der Zypen}
\address{Swiss Armed Forces, CH-3003 Bern,
Switzerland}
\email{dominic.zypen@gmail.com}


\begin{abstract}
	In \cite{Ah}, Ron Aharoni and Vladimir Korman conjectured that
	any hypergraph with only finite edges has a strongly
	minimal cover. We present a counterexample.
\end{abstract}
\parindent = 0mm
\parskip = 2 mm
\maketitle

\section{Introduction}

A {\em hypergraph} is a pair $H=(V,E)$
where $V$ is a set and $E\subseteq \Pow(V)$. The elements of
$E$ are called {\em edges}.

An {\em (edge) cover} of a hypergraph is a subset of 
$K\subseteq E$ such that $\bigcup K = V$. 

We say that a cover $M \subseteq E$ is {\em strongly minimal} if for every
cover $K\subseteq E$ we have $\card(M\setminus K) \leq \card(K\setminus M)$.

\begin{observation} \label{minobs} 
If $M$ is a strongly minimal cover, then it is minimal
with respect to $\subseteq$.\end{observation}
\proof Otherwise suppose that $M_0 \subseteq M$ is a proper subset and
$M_0$ covers the hypergraph. Then 
$\card(M\setminus M_0) > \card(M_0\setminus M) = 0$, so $M$ cannot be
strongly minimal. \proofend

\section{The counterexample}
Conjecture 5.4 of \cite{Ah} says that 
any hypergraph with only finite edges has a strongly
minimal cover. 
\begin{lemma} \label{mainlemma}
	Let $H= (\omega, E)$ where $E$ is the collection of 
	finite subsets of $\omega$. Then $H$ has no strongly minimal
	cover.
\end{lemma}
{\em Proof.} Let $K\subseteq E$ be any cover of $\omega$. We want to 
show that $K$ is not strongly minimal. If $K$ is not minimal, then
by observation \ref{minobs} it cannot be strongly minimal. 

So let us assume
that $K$ is minimal. Recall that $K$ consists only of finite
subsets of $\omega$. As $\bigcup K = \omega$, we can pick
two different members $a, b\in K$. Since $K$ is minimal, we have
$$a\not\subseteq b \text{ and } b\not\subseteq a,$$ and also
$$(a\cup b)\notin K.$$
Let $$K^* = \big(K \setminus \{a, b\}\big) \cup \{a \cup b\}.$$
Clearly $K^*$ is a cover of $\omega$. So
$K \setminus K^* = \{a, b\}$.
On the other hand, $K^*\setminus K$ only has $1$ element, namely $a\cup b$. 

So we have $$\card(K\setminus K^*) = 2 > 1 = \card(K^* \setminus K).$$
Therefore, $K$ is not strongly minimal. \proofend
\section{Acknowledgements}
I want to thank Jonathan David Farley for fruitful discussions
on this topic.
\end{document}